\documentclass[reqno,11pt]{amsart}
\makeatletter
\let\origsection=\section 
\def\section{\@ifstar{\origsection*}{\mysection}}
\def\mysection{\@startsection{section}{1}\z@{.7\linespacing\@plus\linespacing}{.5\linespacing}{\normalfont\scshape\centering\S}}
\makeatother

\usepackage{amsmath,amssymb,amsthm}
\usepackage{mathrsfs}
\usepackage{dsfont}
\usepackage{mathabx}\changenotsign
\usepackage{xparse}
\usepackage{microtype}
\usepackage{xfrac}

\usepackage{xcolor}
\usepackage[backref=page]{hyperref}
\hypersetup{%
    colorlinks,
    linkcolor={red!60!black},
    citecolor={green!60!black},
    urlcolor={blue!60!black}
}

\usepackage{varioref}

\usepackage[english,algoruled,vlined,resetcount,linesnumbered]{
algorithm2e}

\usepackage{bookmark}

\usepackage[abbrev,msc-links,backrefs]{amsrefs}
\usepackage{doi}

\renewcommand{\PrintDOI}[1]{\doi{#1}}

\usepackage[T1]{fontenc}
\usepackage{lmodern}

\usepackage[english]{babel}

\usepackage{tikz}

\usepackage{geometry}
\usepackage{enumitem}

\let\polishlcross=\l
\def\l{\ifmmode\ell\else\polishlcross\fi}

\renewcommand{\setminus}{\smallsetminus}
\renewcommand{\backslash}{\smallsetminus}

\makeatletter
\def\moverlay{\mathpalette\mov@rlay}
\def\mov@rlay#1#2{\leavevmode\vtop{%
   \baselineskip\z@skip \lineskiplimit-\maxdimen
   \ialign{\hfil$\m@th#1##$\hfil\cr#2\crcr}}}
\newcommand{\charfusion}[3][\mathord]{
    #1{\ifx#1\mathop\vphantom{#2}\fi
        \mathpalette\mov@rlay{#2\cr#3}
      }
    \ifx#1\mathop\expandafter\displaylimits\fi}
\makeatother

\newtheoremstyle{case}{}{}{\normalfont}{}{\itshape}{:}{ }{}

\newtheorem{thm}[equation]{Theorem}

\newtheorem{cor}[equation]{Corollary}

\newtheorem{ques}[equation]{Question}

\theoremstyle{definition}

\newtheoremstyle{case}{}{}{\normalfont}{}{\itshape}{\normalfont:}{ }{}

\theoremstyle{case}

\numberwithin{equation}{section}

\let\epsilon\varepsilon
\let\:\colon
\let\subset\subseteq

\def\({\left(}
\def\){\right)}
\def\[{\left[}
\def\]{\right]}
\let\<\langle
\let\>\rangle

\usepackage{datetime}
\usepackage{lineno}
\newcommand*\patchAmsMathEnvironmentForLineno[1]{%
\expandafter\let\csname old#1\expandafter\endcsname\csname #1\endcsname
\expandafter\let\csname oldend#1\expandafter\endcsname\csname end#1\endcsname
\renewenvironment{#1}%
{\linenomath\csname old#1\endcsname}%
{\csname oldend#1\endcsname\endlinenomath}}%
\newcommand*\patchBothAmsMathEnvironmentsForLineno[1]{%
\patchAmsMathEnvironmentForLineno{#1}%
\patchAmsMathEnvironmentForLineno{#1*}}%
\AtBeginDocument{%
\patchBothAmsMathEnvironmentsForLineno{equation}%
\patchBothAmsMathEnvironmentsForLineno{align}%
\patchBothAmsMathEnvironmentsForLineno{flalign}%
\patchBothAmsMathEnvironmentsForLineno{alignat}%
\patchBothAmsMathEnvironmentsForLineno{gather}%
\patchBothAmsMathEnvironmentsForLineno{multline}%
}

\begin{document}

\title{A natural generalisation in Graph Ramsey theory}
\author{
Alexander Haupt
\and Damian Reding
}

\shortdate
\yyyymmdddate
\settimeformat{ampmtime}
\date{\today, \currenttime}

\address{Technische Universit\"at Hamburg, Institut f\"ur Mathematik, Hamburg, Germany}
\email{\{alexander.haupt|damian.reding\}@tuhh.de}

\maketitle

\begin{abstract}
In this note we study graphs $G_r$ with the property that every colouring of $E(G_r)$ with $r+1$ colours admits a copy of some graph $H$ using at most $r$ colours. For $1\le r\le e(H)$ such graphs occur naturally at intermediate steps in the synthesis of a $2$-colour Ramsey graph $G_1\longrightarrow H$. (The corresponding notion of Ramsey-type numbers was introduced by Erd\"os, Hajnal and Rado in 1965 and subsequently studied by Erd\"os and Szemer\'edi in 1972).

For $H=K_n$ we prove a result on building a $G_{r}$ from a $G_{r+1}$ and establish Ramsey-infiniteness. From the structural point of view, we characterise the class of the minimal $G_r$ in the case when $H$ is relaxed to be the graph property of containing a cycle; we then use it to progress towards a constructive description of that class by proving both a reduction and an extension theorem.
\end{abstract}

\vspace{0.5cm}

\section{Introduction and Results}
\label{sec:intro}

The classical $n$-th Ramsey-number $R(n)$ is defined to be the least $N$ such that every $2$-edge-colouring of $G$ admits a copy of $K_n$ using only one of the colours. As opposed to the usual generalisation to $r$-colour-Ramsey numbers, we define the \emph{$n$-th $\sfrac{1}{r}$-Ramsey number $R_{\sfrac{1}{r}}(n)$} to be the least $N$ such that every $(r+1)$-edge-colouring of $G$ admits a copy of $K_n$ using at most $r$ colours. By generalising the classical arguments of~\cite{EG} and~\cite{ES} for $R_{\sfrac{1}{1}}(n)=R(n)$ one obtains

\begin{thm}
Let $r\geq 2$. Then for all $n\geq 2$, $R_{\sfrac{1}{r}}(n)\geq\left[e^{-1}\left(1+\frac{1}{r}\right)+o_r(1)\right] n\left(1+\frac{1}{r}\right)^{\frac{n}{2}}$.
\end{thm}

\begin{thm}
Let $r\geq 2$. Then for all $n\geq 2$, $R_{\sfrac{1}{r}}(n) \le\frac{r(r+2)}{3r+2}\left(1+\frac{1}{r}\right)^{(r+1)n}-r$.
\end{thm}

We supply the arguments in the proofs section. With a little more work (that is essentially different from the $r=1$ case) the upper bound can be improved to $R_{\sfrac{1}{r}}(n) \le r^{\frac{c}{r}n}$, where $c$ is a universal constant (\cite{ES2}). Previously, the quantities had been introduced in~\cite{EHR}.

Our aim in this work is to study the corresponding notion of Ramsey graphs: Given a graph $H$ and $r\in\{1, \ldots, e(H)\}$, we refer to $G$ as a \emph{$\sfrac{1}{r}$-Ramsey graph} for $H$, and write $G\stackrel{\sfrac{1}{r}}{\longrightarrow} H$, if every $(r+1)$-edge-colouring of $G$ admits a copy of $H$ using at most $r$ colours. It is clear that an $\sfrac{1}{r}$-Ramsey graph is $\sfrac{1}{r+1}$-Ramsey and such graphs occur naturally in the sequence
$$H=G_{e(H)}\subset G_{e(H)-1}\subset\ldots\subset G_{2}\subset G_{1}=G,$$
which essentially reduces the construction of $G$ from $H$ to that of $G_r$ from $G_{r+1}$. Note that a necessary condition for $\sfrac{1}{r}$-Ramseyness is the containment of many copies of $H$ that jointly use a relatively low number of edges (indeed, applying the set version of polygamous Hall shows that every $G\stackrel{\sfrac{1}{r}}{\longrightarrow} H$ admits a subgraph $G_0$ containing $>\frac{1}{r} e(G_0)$ copies of $H$).
A construction of $G_{r+1}$ from $G_r$ should therefore seek to maintain this property and we follow up on this observation in the case $H=K_n$. Note that for all $n\geq 3$ and $r\geq 2$, if $G\stackrel{\sfrac{1}{r+1}}{\longrightarrow} K_{R_{\sfrac{1}{r}}(n)}$, then clearly $G\stackrel{\sfrac{1}{r}}{\longrightarrow} K_n$. Indeed, we can reduce the arrowed clique by increasing the number of cliques in $G$ (replace $G$ by the graph $G^+$ obtained from $G$ by adding a new vertex and joining it to every vertex of $G$).
\begin{thm}
For all $n\geq 3$ and $r\geq 2$, if $G\stackrel{\sfrac{1}{r+1}}{\longrightarrow} K_{R_{\sfrac{1}{r}}(n-1, n, \ldots, n)+1}$, then $G^+\stackrel{\sfrac{1}{r}}{\longrightarrow} K_n$.
\end{thm}

The following extends the corresponding result of \cite{NR} for $r=1$.

\begin{thm}
For each $r\geq 2$ there exists $n_r$ such that for all $n\geq n_r$ there exist infinitely many pairwise non-isomorhic minimal $\sfrac{1}{r}$-Ramsey graphs for $K_n$. What's more, we can take $n_2=3$.
\end{thm}

We now relax the basic definition in order to be able to study the structure of $\sfrac{1}{r}$-Ramsey graphs: Given $r\geq 2$, we say $G$ is \emph{$\sfrac{1}{r}$-Ramsey for cyclicity $\mathcal{C}$} (or for any other antimonotone graph property) if every $(r+1)$-colouring of $E(G)$ admits a cycle (of arbitrary length) using at most $r$ colours. We write $\mathcal{R}_{\sfrac{1}{r}}(\mathcal{C})$ for the class of such graphs and $\mathcal{M}_{\sfrac{1}{r}}(\mathcal{C})\subset\mathcal{R}_{\sfrac{1}{r}}(\mathcal{C})$ for the subclass of those graphs, which are minimal wrt. the subgraph relation. This is the terminology adopted from~\cite{RT}, which covers the case $r=1$. Note that as a byproduct of the definition, for $3 \le k \le r$, we have that $C_k\in\mathcal{M}_{\sfrac{1}{r}}(\mathcal{C})$, and also that the diamond graph $K_4-e\in\mathcal{M}_{\sfrac{1}{2}}(\mathcal{C})$. Interestingly, based on the Nash-Williams theorem for multigraphs~\cite{NW, D} a complete characterisation is possible:

\begin{thm}
Let $r\geq 2$. Then for every graph $G$ we have that (a) $G\in\mathcal{R}_{\sfrac{1}{r}} (\mathcal{C})$ if and only if there exists a subgraph $H\subseteq G$ satisfying $re(H)\geq (r+1)v(H)-r$, and (b) $G\in\mathcal{M}_{\sfrac{1}{r}} (\mathcal{C})$ if and only if $re(G)\geq (r+1)v(G)-r$ and $re(H)< (r+1)v(H)-r$ for every proper subgraph $H$ of $G$.
\end{thm}

\begin{cor}
Let $r\geq 2$. Then for every $G\in\mathcal{M}_{\sfrac{1}{r}}(\mathcal{C})$:\\
(a) $e(G)=\left\lceil\left(1+\frac{1}{r}\right)v(G)-1\right\rceil=\left\lfloor\left(1+\frac{1}{r}\right)v(G)-\frac{1}{r}\right\rfloor$.\\
(b) $\delta (G)=2$\\
(c) $v(G)\nequiv 1$ (mod $r$).
\end{cor}

In particular, note that (c) implies that $v(G)$ is even for all $G\in\mathcal{M}_{\sfrac{1}{2}}(\mathcal{C})$.

\begin{thm}
For $r\geq 2$, let $G\in\mathcal{M}_{\sfrac{1}{r}}(\mathcal{C})$ such that $G$ is not a cycle and not the diamond, and let $G'$ denote a graph obtained from $G$ by contracting a shortest cycle in $G$. Then $G'\in\mathcal{R}_{\sfrac{1}{r}}(\mathcal{C})$.
\end{thm}

Note that by alternatingly applying the previous theorem to $G\in\mathcal{M}_{\sfrac{1}{2}}(\mathcal{C})$ and taking minimal $\sfrac{1}{r}$-Ramsey subgraphs, we obtain a proof of:

\begin{cor}
Every graph $G$ satisfying $e(G)\geq\frac{3}{2}v(G)-1$ contains $K_4-e$ as a minor.
\end{cor}

How do we obtain new graphs in $\mathcal{M}_{\sfrac{1}{r}}(\mathcal{C})$ from old ones?

\begin{thm}
For $r\geq 2$, let $G\in\mathcal{M}_{\sfrac{1}{r}}(\mathcal{C})$ be such that $r\notdivides v(G)$, and let $G^+$ denote a graph obtained by a subdividing an arbitrary edge of $G$. Then $G^+\in\mathcal{M}_{\sfrac{1}{r}}(\mathcal{C})$.
\end{thm}

The statement is false if $r\divides v(G)$ (just consider $K_4-e$). To deal with the remaining case, we first note that graphs $G\in\mathcal{M}_{\sfrac{1}{r}}(\mathcal{C})$ are $2$-connected. The following result guarantees that maintaining this condition is already enough to reverse the statement of Theorem 1.7.

\begin{thm}
For $r\geq 2$, let $G\in\mathcal{M}_{\sfrac{1}{r}}(\mathcal{C})$ be such that $r\divides v(G)$, and let $G^+$ denote a $2$-connected graph obtained from $G$ by blowing up an arbitrary vertex to an induced cycle $C$ of length $r+1$. Then $G^+\in\mathcal{M}_{\sfrac{1}{r}}(\mathcal{C})$.
\end{thm}

By iterating the previous theorem we now obtain a constructive proof of the following

\begin{cor}
Let $r\geq 2$ and let $n$ be a proper multiple of $r$. Then there exist $G\in\mathcal{M}_{\sfrac{1}{r}}(\mathcal{C})$ with $v(G)=n$, and in fact such $G$ can be chosen planar. In particular, $\mathcal{M}_{\sfrac{1}{r}}(\mathcal{C})$ contains infinitely many pairwise non-isomorphic planar graphs.
\end{cor}

\section{Proofs}
\label{sec:intro}

\textsl{Proof of Theorem 1.1.}

Colour the edges of $K_N$ uniformly at random with $r+1$ colours. Note that the probability of a fixed copy $K_n\subset K_N$ using at $\leq r$ colours is
$$(r+1)^{-e(K_n)}\sum_{k=1}^r k!\binom{r+1}{k}S(e(K_n); k)\leq(r+1)^{-e(K_n)}\sum_{k=1}^r \binom{r+1}{k} k^{e(K_n)}\leq 2(2^r-1)\left(1+\frac{1}{r}\right)^{-e(K_n)}$$
where $S(e(K_n); k)$ denote the Stirling numbers of the second kind (counting the partitions of $E(K_n)$ into $k$ colour classes). As in Erd\H{o}s' probabilistic bound, we conclude that $R_{\sfrac{1}{r}}(n)>N$, provided
$$\binom{N}{n}\cdot 2(2^r-1)\left(1+\frac{1}{r}\right)^{-\binom{n}{2}}<1\indent\stackrel{\binom{N}{n}\leq\frac{N^n}{n!}}{\Leftarrow}\indent N\geq\left(\frac{n!}{2(2^r-1)}\right)^{\frac{1}{n}}\left(1+\frac{1}{r}\right)^{\frac{1}{2}(n+1)}.$$
Finally, by applying the Stirling bound $n!\leq en^{n+\frac{1}{2}}e^{-n}$, the result follows.\\

\textsl{Proof of Theorem 1.2.}

We will need the auxiliary definition of the corresponding off-diagonal numbers: Define $R_{\sfrac{1}{r}}(n_1,\ldots,n_{r+1})$ to be the least $N$ such that every $(r+1)$-colouring of $E(K_N)$ admits some $K_{n_i}$ missing colour $i$. Note that if some $n_i=2$, then clearly $R_{\sfrac{1}{r}}(n_1,\ldots,n_{r+1})=\min_{j \ne i} \{n_j\}$.\\

\textsl{Claim.} For all $n_1, \ldots, n_{r+1}\geq 3$:\; $R_{\sfrac{1}{r}}(n_1,\ldots,n_{r+1}) \le \left \lceil \frac{1}{r} \sum_{i=1}^{r+1} R_{\sfrac{1}{r}}(\ldots,n_{i-1},n_i-1,n_{i+1},\ldots) \right \rceil$.\\

\textsl{Proof of claim:} Define $R_i:= R_{\sfrac{1}{r}}(\ldots,n_{i-1},n_i-1,n_{i+1},\ldots)$ and $R:= \left\lceil \frac{1}{r} \sum_{i=1}^{r+1} R_i \right \rceil$.

Fix an $(r+1)$-edge-colouring of $K_R$. Fix vertex $x \in K_R$. Suppose that $\forall i \in [r+1]$ at most $R_i-1$ edges incident to $x$ are not using colour $i$; then $\forall i$, at least $(R-1)-(R_i-1)=R-R_i$ are using colour $i$, so we obtain the contradiction:
$$R-1\ge \sum_i (R-R_i) =(r+1)R-\sum_i R_i \ge(r+1)R-rR=R,$$
where the first inequlaity holds by interpreting $R-1$ as the sum of the numbers of the $i$-coloured edges incident to $x$ and the second is due to the ceiling function. Hence, $\exists i\in[r+1]$ such that $\ge R_i$ edges incident to $x$ are not using colour $i$, so there is either $K_{(n_i-1)+1}=K_{n_i}$ not using colour $i$ or some $K_{n_j}$ not using colour $j$ for some $j\in [r+1]\setminus \{i\}$.\\

\textsl{Claim.} Write $N := n_1 + \ldots + n_{r+1}$. There exists $c_r > 0$ such that for all $n_1,\ldots,n_{r+1} \ge 2$,
$$R_{\sfrac{1}{r}}(n_1,\ldots,n_{r+1}) \le c_r \left (1+\frac{1}{r}\right )^N - r$$

\textsl{Proof:}
If all $n_i \ge 3$, suppose (by induction on $N$) that $R_{\sfrac{1}{r}}(\ldots,n_{i-1},n_i-1,n_{i+1},\ldots) \le r_{N-1}$,

$$R_{\sfrac{1}{r}}(n_1,\ldots,n_r)\le \left \lceil \frac{1}{r} \sum_{i=1}^{r+1} r_{N-1} \right \rceil =\left\lceil \left(1+\frac{1}{r}\right) r_{N-1}\right \rceil =\le \left (1+\frac{1}{r}\right) r_{N-1}+1 =: r_N$$

Solving the recursion gives $r_N=c_r \left(1+\frac{1}{r}\right )^N-r$.

Now choose $c_r$ such that the induction start works. Define $M:= \min_{j \ne i} \{n_j\}$. If some $n_i=2$:
$$c_r \left(1+\frac{1}{r}\right)^N-r\ge c_r \left ( 1+\frac{N}{r} \right)-r\ge c_r \left(1+\frac{rM+2}{r} \right )-r=c_r \left (1+\frac{2}{r}+M\right )-r\ge M$$
$$\implies c_r \ge \max_{M \ge 2} \left[ 1+\frac{r^2-r-2}{rM+r+2} \right]=1+\frac{r^2-r-2}{3r+2}=\frac{r(r+2)}{3r+2}$$

\textsl{Proof of Theorem 1.3.}

Fix a colouring $c$: $E(G^+)\longrightarrow [r+1]$. We define an auxiliary colouring $c'$: $E(G)\longrightarrow [r+2]$ as follows: Let $v$ be the new vertex of $G^+$. Give $xy\in E(G)$ colour $i\in [r+1]$ if both $xv, yv$ have colour $i$ in $c$, and give it colour $r+2$ if $xv$, $yv$ have different colours in $c$. By assumption $c'$ now admits a copy of $K:=K_{R_{\sfrac{1}{r}}(n-1, n, \ldots, n)+1}\subset G$ using at most $r+1$ colours. We claim that there exists some colour $i\in [r+1]$ that in $c$ occurs at most once on an edge $vz$ with $z\in V(K)$: otherwise $c'$ admits edges of $K$ of every colour $i\in [r+1]$ and also of colour $r+2$ (namely those joining vertices $z_1, z_2\in V(K)$ with $c(z_1 x)\neq c(z_2 x)$), which is impossible. Hence there exists $K':=K_{R_{\sfrac{1}{r}}(n-1, n, \ldots, n)}\subset K$ with the property that no edge $zv$ with $z\in V(K')$ is using colour $i_0\in [r+1]$ in $c$. Since $R_{\sfrac{1}{r}}(n-1, n, \ldots, n)=R_{\sfrac{1}{r}}(n,\ldots , n-1,\ldots n)$ (with the $n-1$ in $i_0$-th place) there now exists either $j\in [r+1]\backslash\{i_0\}$ with some $K_{n}\subset K$ not using colour $j$ in $c$, or there exists $K_{n-1}\subset K\subset G$, and hence a $K_n$ in $G^+$, none of which is using colour $i_0$.\\

\textsl{Proof of Theorem 1.4.}

Given $r\geq 2$, pick $n_r$ such that $R_{\sfrac{1}{r}}(n)\geq 2n$ for all $n\geq n_r$ (this is possible since $R_{\sfrac{1}{r}}(n)$ grows exponentially). Fix an $n\geq n_r$. Let $F\in\mathcal{M}_{\sfrac{1}{r}}(K_n)$. We find an $H\in\mathcal{M}_{\sfrac{1}{r}}(K_n)$ with $\left|H\right|>\left|F\right|$ by adapting the argument from~\cite{NR}. Let $G\longrightarrow K_n$ with $\chi(H)\leq 2n-1$ for every $H\subset G$ with $\left| H\right|\leq a$, where $a:=\left|F\right|$ (the existence of such $G$ is proven in~\cite{NR}). Clearly $G\in\mathcal{R}_{\sfrac{1}{r}}(K_n)$ and we pick $H\subset G$ minimal with this property. Finally, unless $\left|H\right|>a$, we obtain the contradiction $\chi (H)\geq R_{\sfrac{1}{r}}(n)\geq 2n$ (note that $\chi (H)\geq R_{\sfrac{1}{r}}(n)$ for $H\in\mathcal{R}_{\sfrac{1}{r}}(K_n)$ is obtained by the same simple well-known argument showing that $\chi (H)\geq R(n)$ whenever $H\longrightarrow K_n$).\\
For $r=2$ consider the family of graphs obtained by joining every vertex of an odd cycle to both ends of an independent edge. It is easy to see that each of these graphs is in $\mathcal{R}_{\sfrac{1}{r}}(K_n)$, in fact it is minimal since removing an edge makes the graph $4$-colourable, which then contradicts $\chi (G)\geq R_{\sfrac{1}{2}}(3)=5$. For $n\geq 4$, it is easy to find a colouring of $K_7$ to show that $R_{\sfrac{1}{2}}(4)\geq 8$, whence $R_{\sfrac{1}{2}}(n)\geq 2n$ follows inductively from $R_{\sfrac{1}{2}}(n+1)\geq R_{\sfrac{1}{2}}(n)+2$ (this is proved by adapting the argument for $R(n+1)\geq R(n)+2$). Now the proof of the first part applies.\\

\newpage
\textsl{Proof of Theorem 1.5.}

Note that, since the statement in (b) characterizes those Ramsey graphs with no proper Ramsey subgraphs, it follows immediately from (a), which we now prove. For the \emph{if} direction, suppose that $G$ is not a $\sfrac{1}{r}$-Ramsey graph for $\mathcal{C}$, so there exists a $(r+1)$-colouring of the edges of $G$ with the property that the union of any $r$ of the $r+1$ colour classes $E_1,\ldots, E_{r+1}$ is an acyclic graph on $V(G)$. Towards a contradiction, we can thus write
$$re(H) = r\sum\limits_{i=1}^{r+1}\left|E(H)\cap E_i\right| = \sum\limits_{j=1}^{r+1}\sum\limits_{i\neq j}\left|E(H)\cap E_i\right|\leq (r+1)(v(H)-1) < (r+1)v(H)-r$$
For the \emph{only if} direction, suppose that $re(H)\leq (r+1)v(H)-r-1=(r+1)(v(H)-1)$ for every subgraph $H\subseteq G$; this condition enables us to apply the Nash-Williams theorem~\cite{NW} to the multigraph $G^r$ obtained from $G$ by replacing each edge of $G$ by $r$ edges. We obtain a decomposition $E(G)=E_1\cup\ldots\cup E_r$ into simple forests on $V(G)$, so no two of the $r$ edges joining a pair of adjacent vertices $x, y$ in $G^r$ are of the same colour. Hence each multiedge $e^r\in E(G^r)$ must be using precisely $r$ of the $r+1$ colours, so we may define an edge-colouring of $G$ by choosing for the corresponding edge $e\in E(G)$ the unique missing colour each. In this colouring, given the acyclicity of the $E_i$, every cycle in $G$ must be using all $r+1$ colours.\\

\textsl{Proof of Corollary 1.6.}

The previous theorem yields $re(G)\geq (r+1)v(G)-r$ and $re(G-e)< (r+1)v(G-e)-r$, or equivalently $r(e(G)-1)\leq(r+1)(v(G)-1)$. Solving this for $e(G)$ gives

$$\left\lceil\left(1+\frac{1}{r}\right)v(G)-1\right\rceil\leq e(G)\leq\left\lfloor\left(1+\frac{1}{r}\right)v(G)-\frac{1}{r}\right\rfloor$$

with both bounds in fact equal (write $v(G)=kr+d$, $0\leq d<r$). The average degree of $G$ is

$$\frac{2e(G)}{v(G)}\leq\frac{2}{v(G)}\left(\left(1+\frac{1}{r}\right)v(G)-\frac{1}{r}\right)= 2+\frac{2}{r}\left(1-\frac{1}{v(G)}\right)\leq 3-\frac{1}{v(G)}<3,$$

hence there exists a vertex of degree $\leq 2$. The second result follows since clearly $\delta(G)\geq 2$. Finally, suppose that $v(G)=1+dr$, $d\in\mathbb{N}$, for some $G\in\mathcal{M}_{\sfrac{1}{r}}(\mathcal{C})$. Then $e(G)=dr+d+1$ by (a). By (b) there exists $v\in V(G)$; consider $H:=G-v$. Then $re(H)\leq (r+1)v(H)-r-1$, where $e(H)=dr+d-1$ and $v(H)=dr$. This results in the contradiction $0\leq -1$.\\

\textsl{Proof of Theorem 1.7.}

Since a shortest cycle in $G$ is necessary induced, we have that $v(G')=v(G)-(g-1)$ and $e(G')=e(G)-g$, where $g$ denotes the girth of $G$. Note that $g\geq r+1$ by minimality.
\begin{eqnarray*}
re(G') & = & r(e(G)-g)\\
& \geq & (r+1)v(G)-r-rg\\
& = & (r+1)(v(G')+(g-1))-r-rg\\
& = & (r+1)v(G')-r+(g-r-1)\\
& \geq & (r+1)v(G')-r
\end{eqnarray*}

Note that while the argument itself does not require that $G\neq D_r$, we disallow this graph on the understanding that $D_r'$ is then actually a multigraph (the double edge between two vertices) and hence not a member of $\mathcal{M}_{\sfrac{1}{r}}(\mathcal{C})$.\\

\textsl{Proof of Theorem 1.9.}

To begin with, first note that graphs $G\in\mathcal{M}_{\sfrac{1}{r}}(\mathcal{C})$ are $2$-connected (indeed, if there is a cutvertex $v\in V(G)$, let $C$ be a component of $G-v$; put an $(r+1)$-edge-colouring on $G[V(C)\cup\{v\}]$ and $G[V\backslash V(C)]$ without cycles using $\leq r$ colours. Then, in the resulting edge-colouring of $G$, all cycles use all the colours).\\

Let $\ell \in [r-1]$ with $v(G) \equiv \ell$ (mod $r$).
We have $e(G^+)=e(G)+1$ and $v(G^+)=v(G)+1$, so
\begin{eqnarray*}
	G\in\mathcal{M}_{\sfrac{1}{r}}(\mathcal{C}) &\implies& r e(G) = r \cdot \left\lceil\left(1+\frac{1}{r}\right)v(G)-1\right\rceil = (r+1)v(G)-\ell \\
	&\implies& r (e(G)+1) = (r+1)(v(G)+1)-\ell-1\\
	&\implies& r e(G^+) = (r+1)v(G^+)-\ell-1\\
	&\implies& re(G^+) \ge (r+1)v(G^+)-(r-1)-1\\
	&\implies& re(G^+) \ge (r+1)v(G^+)-r
\end{eqnarray*}

Let $e_0 \in E(G) \setminus E(G^+)$ be the edge that was subdivided and $e_1,e_2 \in E(G^+) \setminus E(G)$ be the new edges.

To show minimality, it suffices to show that $\forall e \in E(G^+)$ we have $G^+ - e \notin \mathcal{R}_{\sfrac{1}{r}}(\mathcal{C})$. If $e=e_1$ or $e=e_2$, let $c$ be a $(r+1)$-colouring of $G-e_0$ and set $c(e_2)=1$ or $c(e_1)=1$ respectively. Finally, if $e \in E(G-e_0)$, let $c$ be a $\sfrac{1}{r}$-edge-colouring of $G-e$. Define $c' : E(G^+-e) \to [r+1]$ by $c'(e)=c(e)$ if $e\in E(G)$ and $c'(e)=c(e_0)$ otherwise. In all cases we have found a $\sfrac{1}{r}$-edge-colouring of $G^+-e$ and thus $G^+ \in \mathcal{M}_{\sfrac{1}{r}}(\mathcal{C})$.\\

\textsl{Proof of Theorem 1.10.}

We have $e(G^+)=e(G)+r+1$ and $v(G^+)=v(G)+r$, so
\begin{eqnarray*}
G\in\mathcal{M}_{\sfrac{1}{r}}(\mathcal{C}) &\implies& re(G) \ge (r+1)v(G)-r \\
&\implies& r(e(G)+(r+1)) \ge (r+1)(v(G)+r)-r \\
&\implies& re(G^+) \ge (r+1)v(G^+)-r
\end{eqnarray*}
Now let $H^+$ be a proper subgraph of $G^+$ and let $H$ denote the subgraph of $G$ consisting of those vertices and edges, which belong to $H^+$ prior to the contraction of $C$. Set $k := e(H^+)-e(H)$. Clearly $0\le k \le r+1$ and $v(H^+)-v(H) \ge k$.\\
If $H$ is a proper subgraph of $G$, then we have that $v(H^+)-v(H) \ge rk/(r+1)$, so
$$re(H^+)=re(H)+rk < (r+1)v(H)-r+rk  \le (r+1)v(H^+)-r$$
If $H=G$, we further distinguish between the cases $k=0$ and $1 \le k \le r$:\\
\indent If $k=0$, we have that $v(H^+)-v(G) = \left | \{v\in C : d(v)\ge 3\} \right | - 1 \ge 1$ by $2$-connectedness.\\
\indent If $1 \le k \le r$, we have $v(H^+)-v(G) \ge k$. Note that $k = r+1$ cannot happen, as $H^+$ is a proper subgraph of $G^+$.

Thus, for all values of $k$, we have $(r+1) (v(H^+)-v(G)) > rk$.

It follows that $re(H^+)=re(G)+rk=(r+1)v(G)-r+rk < (r+1)v(H^+)-r$.

Therefore $re(H^+) < (r+1)v(H^+)-r$ and so $G^+\in\mathcal{M}_{\sfrac{1}{r}}(\mathcal{C})$.

\section{Concluding Remarks}
\label{sec:intro}

It is clear that for $r\geq\binom{n}{2}$ the only minimal $\sfrac{1}{r}$-Ramsey graph for $K_n$ is $K_n$ itself. Theorem 1.4 guarantees that for $n\geq n_r$, however, $K_n$ is Ramsey-infinite and an upper bound on $n_r$ could be obtained by tracking the computations. But what happens when $r$ is large relative to $n$?

\begin{ques}
Is $K_n$ $\sfrac{1}{r}$-Ramsey-infinite for all $r<\binom{n}{2}$?
\end{ques}

Further, towards the purpose of building classical Ramsey graphs from $\sfrac{1}{r}$-Ramsey graphs it would be interesting to determine the minimum size of a clique that a graph needs to $\sfrac{1}{r+1}$-arrow in order for it to $\sfrac{1}{r}$-arrow a smaller clique of given size: We define $f_{r; k}(n)$ to be the minimum $N$ such that for every graph $G\in\mathcal{R}_{\sfrac{1}{r+1}}(K_{f_{r; k}(n)})$ we have that $G^{+k}\in\mathcal{R}_{\sfrac{1}{r}}(K_n)$, where $G^{+k}$ is the graph obtained recursively from $G$ by joining a new vertex to every vertex of $G^{+(k-1)}$. Proposition 1.3 shows that $f_{r; 1}(n)\leq R_{\sfrac{1}{r}}(n-1, n,\ldots , n)$.

\begin{ques}
Determine $f_{r; k}(n)$.
\end{ques}

On the subject of cyclicity, observe that our results imply that, for $n\geq 3$, there exist $G\in\mathcal{M}_{\sfrac{1}{2}}(\mathcal{C})$ on $n$ vertices if and only if $n$ is even. For $r\geq 3$ in turn, with regard to the existence of $G\in\mathcal{M}_{\sfrac{1}{2}}(\mathcal{C})$ with $v(G)=n$ we know that it is necessary for $n$ to satisfy $v(G)\nequiv 1$ (mod $r$), and that it is sufficient for $n$ to be one of $3,\ldots, r-1$ (in which case $G=C_n$) or a multiple of $r$, but that in general this sufficient condition itself is not necessary (consider the graph obtained from $K_5$ by removing the $5$-vertex graph $K_3+K_2$). We therefore raise the following

\begin{ques}
Given $r\geq 3$, for which $n\in\mathbb{N}$ exactly does there exist $G\in\mathcal{M}_{\sfrac{1}{r}}(\mathcal{C})$ on $n$ vertices?
\end{ques}

Finally, we remark that by adapting the above methods one could in principle initiate a study of graphs, for any fixed $r<s$, with the property that every colouring of $E(G)$ with $s$ colours admits a copy of some graph $H$ using at most $r$ colours (with $r=1$ being the classical multicolour case). However, in general such graphs seem to exhibit more irregular behaviour than the particular case $s=r+1$ (e.g. for $s=4$ and $r=2$ there exist minimal graphs $G_1$ and $G_2$ with $v(G_1)=v(G_2)$, but $e(G_1)\neq e(G_2)$). It is amusing that the cube graph on $8$ vertices turns out to be a minimal graph for $s=5$ and $r=3$.\\\\
\noindent\emph{Acknowledgements.} We thank Dennis Clemens for pointing out a first proof of Corollary 1.6(c) and Yoshiharu Kohayakawa for pointing out the reference~\cite{ES2}.

\bibliographystyle{abbrv}
\bibliography{partialramsey}

\end{document}